\input pstricks\input xy \xyoption{all}

\newwrite\num
%
\output={\if N\header\headline={\hfill}\fi
\plainoutput\global\let\header=Y}
\magnification\magstep1
\tolerance = 500
\hsize=14.4true cm
\vsize=22.5true cm
\parindent=6true mm\overfullrule=2pt
\newcount\kapnum \kapnum=0
\newcount\parnum \parnum=0
\newcount\procnum \procnum=0
\newcount\nicknum \nicknum=1
\font\saser=cmss10
\font\ninett=cmtt9

\font\ninebf=cmbx9

\font\sixbf=cmbx6
\font\ninesl=cmsl9

\font\nineit=cmti9

\font\ninerm=cmr9

\font\sixrm=cmr6
\font\ninei=cmmi9
\font\eighti=cmmi8
\font\sixi=cmmi6
\skewchar\ninei='177 \skewchar\eighti='177 \skewchar\sixi='177
\font\ninesy=cmsy9
\font\eightsy=cmsy8
\font\sixsy=cmsy6
\skewchar\ninesy='60 \skewchar\eightsy='60 \skewchar\sixsy='60
\font\titelfont=cmr10 scaled 1440
\font\paragratit=cmbx10 scaled 1200

\font\name=cmcsc10
\font\emph=cmbxti10

\font\tenmsbm=msbm10
\font\sevenmsbm=msbm7
%

%
\font\got=eufm10

\font\teneufm=eufm10
\font\seveneufm=eufm7
\font\fiveeufm=eufm5
\newfam\eufmfam
\textfont\eufmfam=\teneufm
\scriptfont\eufmfam=\seveneufm
\scriptscriptfont\eufmfam=\fiveeufm

\font\tenmsam=msam10
\font\sevenmsam=msam7
\font\fivemsam=msam5
\newfam\msamfam
\textfont\msamfam=\tenmsam
\scriptfont\msamfam=\sevenmsam
\scriptscriptfont\msamfam=\fivemsam
\font\tenmsbm=msbm10
\font\sevenmsbm=msbm7
\font\fivemsbm=msbm5
\newfam\msbmfam
\textfont\msbmfam=\tenmsbm
\scriptfont\msbmfam=\sevenmsbm
\scriptscriptfont\msbmfam=\fivemsbm
\def\Bbb#1{{\fam\msbmfam\relax#1}}
\def\cz{{\kern0.4pt\Bbb C\kern0.7pt}
}
\def\ez{{\kern0.4pt\Bbb E\kern0.7pt}
}
\def\fz{{\kern0.4pt\Bbb F\kern0.3pt}}
\def\gz{{\kern0.4pt\Bbb Z\kern0.7pt}}
\def\hz{{\kern0.4pt\Bbb H\kern0.7pt}
}
\def\kz{{\kern0.4pt\Bbb K\kern0.7pt}
}
\def\nz{{\kern0.4pt\Bbb N\kern0.7pt}
}
\def\oz{{\kern0.4pt\Bbb O\kern0.7pt}
}
\def\rz{{\kern0.4pt\Bbb R\kern0.7pt}
}
\def\sz{{\kern0.4pt\Bbb S\kern0.7pt}
}
\def\pz{{\kern0.4pt\Bbb P\kern0.7pt}
}
\def\qz{{\kern0.4pt\Bbb Q\kern0.7pt}
}
\newskip\ttglue
\def\ninepoint{\def\rm{\fam0\ninerm}%
  \textfont0=\ninerm \scriptfont0=\sixrm \scriptscriptfont0=\fiverm
  \textfont1=\ninei \scriptfont1=\sixi \scriptscriptfont1=\fivei
  \textfont2=\ninesy \scriptfont2=\sixsy \scriptscriptfont2=\fivesy
  \textfont3=\tenex \scriptfont3=\tenex \scriptscriptfont3=\tenex
  \def\it{\fam\itfam\nineit}%
  \textfont\itfam=\nineit
  \def\sl{\fam\slfam\ninesl}%
  \textfont\slfam=\ninesl
  \def\bf{\fam\bffam\ninebf}%
  \textfont\bffam=\ninebf \scriptfont\bffam=\sixbf
   \scriptscriptfont\bffam=\fivebf
  \def\tt{\fam\ttfam\ninett}%
  \textfont\ttfam=\ninett
  \tt \ttglue=.5em plus.25em minus.15em
  \normalbaselineskip=11pt
  \font\name=cmcsc9
  \let\sc=\sevenrm
  \let\big=\ninebig
  \setbox\strutbox=\hbox{\vrule height8pt depth3pt width0pt}%
  \normalbaselines\rm
  \def\sl{\it}}

\headline={\ifodd\pageno\rightheadline\else\leftheadline\fi}
\def\rightheadline{\ninepoint Paragraphen"uberschrift\hfill\folio}
\def\leftheadline{\ninepoint\folio\hfill Chapter"uberschrift}
\let\header=Y
\def\titel#1{\need 9cm \vskip 2truecm
\parnum=0\global\advance \kapnum by 1
{\baselineskip=16pt\lineskip=16pt\rightskip0pt
plus4em\spaceskip.3333em\xspaceskip.5em\pretolerance=10000\noindent
\titelfont Chapter \uppercase\expandafter{\romannumeral\kapnum}.
#1\vskip2true cm}\def\leftheadline{\ninepoint
\folio\hfill Chapter \uppercase\expandafter{\romannumeral\kapnum}.
#1}\let\header=N
}
\def\Titel#1{\need 9cm \vskip 2truecm
\global\advance \kapnum by 1
{\baselineskip=16pt\lineskip=16pt\rightskip0pt
plus4em\spaceskip.3333em\xspaceskip.5em\pretolerance=10000\noindent
\titelfont\uppercase\expandafter{\romannumeral\kapnum}.
#1\vskip2true cm}\def\leftheadline{\ninepoint
\folio\hfill\uppercase\expandafter{\romannumeral\kapnum}.
#1}\let\header=N
}
\def\need#1cm {\par\dimen0=\pagetotal\ifdim\dimen0<\vsize
\global\advance\dimen0by#1 true cm
\ifdim\dimen0>\vsize\vfil\eject\noindent\fi\fi}
\def\neupara#1{\par\penalty-2000
\procnum=0\global\advance\parnum by 1
\vskip1cm\noindent{\paragratit \the\parnum. #1}%
\def\rightheadline{\ninepoint\S\the\parnum.\ #1\hfill \folio}%
\vskip 8mm\noindent}
\def\Proclaim #1 #2\finishproclaim {\bigbreak\noindent
{\bf#1\unskip{}. }{\it#2}\medbreak\noindent}
%
\gdef\proclaim #1 #2 #3\finishproclaim {\bigbreak\noindent%
\global\advance\procnum by 1
{%
{\relax\ifodd \nicknum
\hbox to 0pt{\vrule depth 0pt height0pt width\hsize
   \quad \ninett#3\hss}\else {}\fi}%
\bf\the\parnum.\the\procnum\ #1\unskip{}. }
{\it#2}
\immediate\write\num{\string\def
 \expandafter\string\csname#3\endcsname
 {\the\parnum.\the\procnum}}
\medbreak\noindent}
\newcount\stunde \newcount\minute \newcount\hilfsvar
\def\uhrzeit{
    \stunde=\the\time \divide \stunde by 60
    \minute=\the\time
    \hilfsvar=\stunde \multiply \hilfsvar by 60
    \advance \minute by -\hilfsvar
    \ifnum\the\stunde<10
    \ifnum\the\minute<10
    0\the\stunde:0\the\minute~Uhr
    \else
    0\the\stunde:\the\minute~Uhr
    \fi
    \else
    \ifnum\the\minute<10
    \the\stunde:0\the\minute~Uhr
    \else
    \the\stunde:\the\minute~Uhr
    \fi
    \fi
    }

\def\calG{{\cal G}} \def\calH{{\cal H}}

\def\gota{\hbox{\got a}} 
\def\gotb{\hbox{\got b}} 
\def\gotc{\hbox{\got c}}

\def\gotx{\hbox{\got x}}

\def\GL{\mathop{\rm GL}\nolimits}

\def\kernel{\mathop{\rm kernel}\nolimits}

\def\mod{\mathop{\rm mod}\nolimits}

\def\proj{\mathop{\rm proj}\nolimits}

\def\Sp{\mathop{\rm Sp}\nolimits}

\def\boxit#1{
  \vbox{\hrule\hbox{\vrule\kern6pt
  \vbox{\kern8pt#1\kern8pt}\kern6pt\vrule}\hrule}}
\def\Boxit#1{
  \vbox{\hrule\hbox{\vrule\kern2pt
  \vbox{\kern2pt#1\kern2pt}\kern2pt\vrule}\hrule}}

\def\zwischen#1{\bigbreak\noindent{\bf#1\medbreak\noindent}}

\def\smallni{\smallskip\noindent }
\def\medni{\medskip\noindent }
\def\bigni{\bigskip\noindent }

\def\lo{\longrightarrow}

\def\loma{\longmapsto}

\def\spitz#1{\langle#1\rangle}
\def\imag{{\rm i}}
\def\pii{\pi {\rm i}}

\def\set#1{\bigl\{\,#1\,\bigr\}}

\def\square{\hbox{\hbox to 0pt{$\sqcup$\hss}\hbox{$\sqcap$}}}
\def\qed{\ifmmode\square\else{\unskip\nobreak\hfil
\penalty50\hskip3em\null\nobreak\hfil\square
\parfillskip=0pt\finalhyphendemerits=0\endgraf}\fi}
\def\pn{\the\parnum.\the\procnum}
\def\downmapsto{{\buildrel
        {\vbox{\hbox{\hskip.2pt$\scriptstyle-$}}}
        \over{\raise7pt\vbox{\vskip-4pt\hbox{$\textstyle\downarrow$}}}}}

\def\G2Pr{4.1}

\def\Goep{5.2}
\nopagenumbers
\immediate\newwrite\num
\nicknum=0  
\let\header=N

\immediate\openout\num=goepel.num
\immediate\newwrite\num\immediate\openout\num=goepel.num
\def\RAND#1{\vskip0pt\hbox to 0mm{\hss\vtop to 0pt{%
\raggedright\ninepoint\parindent=0pt%
\baselineskip=1pt\hsize=2cm #1\vss}}\noindent}
\noindent
\centerline{\titelfont The G\"opel variety}%
\bigskip
\centerline{2017}
\def\leftheadline{\ninepoint\folio\hfill
The G\"opel variety}%
\def\rightheadline{\ninepoint Introduction\hfill \folio}%
\headline={\ifodd\pageno\rightheadline\else\leftheadline\fi}
\vskip 1.5cm
\leftline{\it \hbox to 6cm{Eberhard Freitag\hss}
Riccardo Salvati
Manni  }
\leftline {\it  \hbox to 6cm{Mathematisches Institut\hss}
Dipartimento di Matematica, }
\leftline {\it  \hbox to 6cm{Im Neuenheimer Feld 288\hss}
Piazzale Aldo Moro, 2}
\leftline {\it  \hbox to 6cm{D69120 Heidelberg\hss}
I-00185 Roma, Italy. }
\leftline {\tt \hbox to 6cm{freitag@mathi.uni-heidelberg.de\hss}
salvati@mat.uniroma1.it}
\vskip1cm\noindent%
\let\header=N%
\def\imag{{\rm i}}%
\def\transpose#1{\kern1pt{^t\kern-3pt#1}}%
\centerline{\vbox{\noindent\hsize=10cm
{\ninepoint{\bf Abstract}
\smallni
In this paper we will  prove that the  six-dimensional {\it G\"opel variety\/} in $P^{134}$
is generated by 120 linear, 35 cubic and 35 quartic relations.
This  result was already  obtained  in [RS] , but the authors  used a statement  in  [Co] saying
that the G\"opel variety set theoretically is generated by the linear and cubic relations alone.
Unfortunately this statement is false. There are 120  extra points. Nevertheless
the results  stated  in   [RS] are correct.  There are required several
changes that  we will illustrate in some detail.}}}
\bigni
{\paragratit Introduction}%
\medni
In this paper we will  consider the so-called, in [RS], {\it G\"opel variety\/}.
It  is an interesting six  dimensional  variety arising in  connection  with
the  moduli space of  principally  polarized abelian  varieties of level 2 and genus 3.
In  fact, these two  varieties are  birational equivalent,
since the G\"opel variety parameterizes the coefficients of the   Coble quartic.
Moreover, this  can  be seen as the   the {\saser GIT} quotient   describing the
moduli space of  seven points in  $P^2$, i.e.\
the moduli space of marked  Del Pezzo surfaces of degree 2 , cf.~[DO] .
It  has  also a realization as  ball quotient, cf.~[Ko].
We start  with a description.
\smallskip
We consider the vector space of ``theta characteristics'' $\fz_2^{2g}$ equipped with the standard
symplectic pairing. A {\it G\"opel group\/} is a maximal isotropic subspace $G$. Its dimension is $g$.
Every G\"opel group has exactly only translate $M=a+G$ that consists of even characteristics.
We ssociate two modular forms (Sect.~1)
$$\vartheta_G=\prod_{m\in M}\vartheta[m],\quad s_G=\prod_{m\notin M}\vartheta[m].$$
These are modular forms of of level two. In the cases $g\ge 3$ the multipliers are trivial.
The graded algebras and their associated projective varieties
$$A=\cz[\dots \vartheta_G\dots],\quad B=\cz[\dots s_G\dots]$$
are of great interest. We are mainly interested in the case $g=3$. The 6-dimensional variety
$\proj(B)$ is   the {\it G\"opel variety}. It has been studied several times.
There are 135 G\"opel groups. Hence this variety sits in $P^{134}$. As a consequence of the
quartic Riemann theta relations one obtains 120 linear relations between the $s_G$. Hence the G\"opel
variety can be considered as a subvariety of $P^{14}$. In the literature there have been described also
cubic [Co], [DO] and quartic relations [RS] and it is stated in this paper
that these linear cubic and
quartic relations generate the full
ideal of relations, cf. [RS], Theorem 5.1. As far as we understand, the proof of this fact depends on a statement of
Coble that the G\"opel variety set theoretically is generated by the linear and cubic relations alone.
Unfortunately this statement is false. There are 120 extra points that satisfy the linear and cubic relations
but not all quartic relations. Nevertheless it remains true that the linear, cubic and quartic relations
generate the full vanishing ideal (hence a prime ideal). The proof is quite involved, so we decided to
reopen the story again and to give a complete proof. Needless to say that nowadays computers have been
a basic tool.  For sake  of completeness  we restate  the part of  Theorem 5.1 from [RS] that  we will  "reprove".
\Proclaim
{Theorem}
{The six-dimensional G\"opel variety $\calG$ has degree 175 in $P^{14}$. The homogeneous
coordinate ring of $\calG$ is Gorenstein
and its defining prime ideal is minimally generated by 35 cubics and 35 quartics.}
\finishproclaim
As in [RS] computer computations are necessary. We used the computer algebra system
{\saser MAGMA} [BC].
\neupara{Thetanullwerte}%
In our context, theta characteristics in genus $g$ are elements of $\fz_2^{2g}$.
Usually they are written as columns.
Sometimes we have to associate to an element of $\fz_2$ an integer. This is done by means of a
section
$$\fz_2\buildrel\iota\over\lo\gz\lo\fz_2$$
We take
$$\iota(x)=\cases{0\in\gz&if $x=0$ in $\fz_2$,\cr1\in\gz&if $x=1$ in $\fz_2$.\cr}$$
A characteristic $m$ is called even if
$$a'b=0\quad\hbox{where}\quad m={a\choose b}.$$
There are $(2^g+1)2^{g-1}$ even characteristics.
The full modular group $\Gamma_g:=\Sp(g,\gz)$ acts on the set of (even) characteristics
by
$$M\{m\}=M'^{-1}m+\pmatrix{(CD')_0\cr (AB')_0}\qquad (\mod\; 2).$$
Here we denote by $S_0$ the diagonal of the matrix $S$ written as column.
This action is double transitive on the set of all even characteristics.
\medskip
Thetanullwerte are defined by
$$\vartheta[m](Z)=\sum_{n\in\gz^n} e^{\pii\{Z[n+a/2]+b'(2n+a)\}},\qquad m={a\choose b}.$$
Here $m$ is a characteristic (to be precise, one has it to replace be its image in $\gz^{2g}$
under the section $\iota$).  A thetanullwert is different from zero if and only if $m$ is
even. Hence we have $(2^g+1)2^{g-1}$ non-zero thetanullwerte.
The thetanullwerte satisfy the {\it quartic Riemann relations.}
They are of the following type. For each pair of characteristics $p,q\in\fz^{2g}$ one has
$$\sum_{m\;\hbox{\sevenrm even}}\pm\vartheta[m]\vartheta[m+p]\vartheta[m+q]\vartheta[m+p+q]=0.$$
We will not describe the quite delicate signs here.
\medskip
We number the 36 even theta characteristics in genus 3 as follows:
\bigskip
\hbox{
\vbox{\halign{\qquad$#$\quad\hfill&$#$\hfill\cr
1
&(0,0,0,0,0,0)\cr
2
&(1,0,0,0,0,0)\cr
3
&(0,1,0,0,0,0)\cr
4
&(1,1,0,0,0,0)\cr
5
&(0,0,1,0,0,0)\cr
6
&(1,0,1,0,0,0)\cr
7
&(0,1,1,0,0,0)\cr
8
&(1,1,1,0,0,0)\cr
9
&(0,0,0,1,0,0)\cr
10
&(0,1,0,1,0,0)\cr
11
&(0,0,1,1,0,0)\cr
12
&(0,1,1,1,0,0)\cr}}\qquad
\vbox{\halign{$#$\quad\hfill&$#$\hfill\cr
13
&(0,0,0,0,1,0)\cr
14
&(1,0,0,0,1,0)\cr
15
&(0,0,1,0,1,0)\cr
16
&(1,0,1,0,1,0)\cr
17
&(0,0,0,1,1,0)\cr
18
&(1,1,0,1,1,0)\cr
19
&(0,0,1,1,1,0)\cr
20
&(1,1,1,1,1,0)\cr
21
&(0,0,0,0,0,1)\cr
22
&(1,0,0,0,0,1)\cr
23
&(0,1,0,0,0,1)\cr
24
&(1,1,0,0,0,1)\cr}}\qquad
\vbox{\halign{$#$\quad\hfill&$#$\hfill\cr
25
&(0,0,0,1,0,1)\cr
26
&(0,1,0,1,0,1)\cr
27
&(1,0,1,1,0,1)\cr
28
&(1,1,1,1,0,1)\cr
29
&(0,0,0,0,1,1)\cr
30
&(1,0,0,0,1,1)\cr
31
&(0,1,1,0,1,1)\cr
32
&(1,1,1,0,1,1)\cr
33
&(0,0,0,1,1,1)\cr
34
&(1,1,0,1,1,1)\cr
35
&(1,0,1,1,1,1)\cr
36
&(0,1,1,1,1,1)\cr
}}}
\bigni
\neupara{Thetanullwerte and modular forms}%
Recall that the real symplectic group
$$\Sp(g,\rz)=\set{M\in\GL(2g,\rz);
\quad {M'}IM=I}\qquad \left(I=\pmatrix{0&-E\cr E&0}\right)$$
acts on the generalized half plane
$$\calH_g:=\set{Z=X+\imag Y;\quad Z={Z'},\  Y>0\ \hbox{(positive definite)}}$$
by
$$MZ=(AZ+B)(CZ+D)^{-1},\qquad M=\pmatrix{A&B\cr C&D}.$$
Let $\Gamma_g:=\Sp(g,\gz)$ be  the Siegel modular group.
The principal congruence subgroup of level $l$ is
$$\Gamma_g[l]:=\kernel (\Sp(g,\gz)\lo\Sp(g,\gz/l\gz))$$
and Igusa's subgroup is
$$\Gamma_g[l,2l]:=\set{M\in\Gamma_g[l];\quad
A{B'}/l\ \hbox{and}\ C{D'}/l\ \hbox{have even diagonal}}.$$
For even $l$, Igusa's subgroup  is a normal subgroup of $\Gamma_g$.
\smallskip
For each $M\in\Sp(g,\rz)$ and each $Z\in\calH_g$ we take a holomorphic
root $\sqrt{\det(CZ+D)}$. To be concrete, we make the choice such that for $Z=\imag E$
it is the principal value of the square root.
Let $\Gamma\subset\Sp(g,\gz)$ be a subgroup of finite index and $r$ an integer. A multiplier system
of weight $r/2$ is map $v:\Gamma\to S^1$ such that
$v(M)\sqrt{\det(CZ+D)}^{\,r/2}$ is an automorphy factor. This means
that
$$(f\vert M)(Z):=v(M)^{-1}\sqrt{\det(CZ+D)}^{\,-r/2}f(MZ),\quad M\in\Gamma,$$
is an action of $\Gamma$ (from the right). A modular form of weight $r/2$ and multiplier system
$v$ on $\Gamma$ is a holomorphic function $f:\calH_g\to\cz$ with the property
$$f(MZ)=v(M)\sqrt{\det(CZ+D)}^{\>r}f(Z),\qquad M\in\Gamma.$$
In the case $g=1$ the usual regularity condition at the cusps has to be added.
We denote the space of all these forms by $[\Gamma,r/2,v]$. In the case that $r/2$ is integral and
that $v$ is trivial we omit $v$ in this notation.
For given $(r_0,v_0)$ we can define the algebra
$$A(\Gamma,(r_0,v_0)):=\bigoplus_{r\in\gz}[\Gamma,rr_0,v^r].$$
If it is clear which starting weight and multiplier system are used, we write
$$A(\Gamma)=A(\Gamma,(r_0,v_0)).$$
By the theory of Satake compactification, $A(\Gamma)$ is a finitely generated algebra
whose associated projective variety, considered as complex space, is biholomorphic
equivalent to the Satake compactification of $\calH_g/\Gamma$,
$$\proj(A(\Gamma))=X(\Gamma)=\overline{\calH_g/\Gamma}.$$
There is a fundamental multiplier system on the theta group
$\Gamma_\vartheta:=\Gamma_g[1,2]$,
of 8-th roots of unity $v(M)$ for $M\in\Gamma_\vartheta$, such that
$$\vartheta[0]\in[\Gamma_\vartheta,1/2,v_\vartheta].$$
More generally, a formula
$$\vartheta[M\{m\}](MZ)=v(M,m)\sqrt{\det(CZ+D)}\vartheta[m](Z)\quad\hbox{for all}\ M\in\Gamma_g,$$
with a certain system of numbers $v(M,m)$ of 8th roots of unity holds.
\smallskip
It is well-known that
$$\vartheta[m]\in[\Gamma_g[4,8],1/2,v_\vartheta]\quad\hbox{for all}\ m.$$
We use $(1/2,v_\vartheta)$ as starting weight and multiplier system to define the ring
$$A(\Gamma_g[4,8]):=\bigoplus_{r\in\gz}[\Gamma_g[4,8],r/2,v_\vartheta^r].$$
The fundamental lemma of Igusa [Ig1] says that this ring is the normalization of the ring generated by the
thetanullwerte. In the cases $g\le 2$ both rings agree [Ig2].
\neupara{Thetanullwerte of second kind}%
The $2^g$ thetanullwerte of second kind are defined by
$$f_a(Z):=\vartheta\Bigl[{a\atop 0}\Bigr](2Z)=\sum_{n\in\gz^g}e^{2\pii Z[n+a/2]}.$$
They are modular forms on the group
$\Gamma_g[2,4]$ with a joint multiplier system $v_\Theta$ on $\Gamma_g[2,4]$. One has
$v_\Theta^2=v_\vartheta^2$ on $\Gamma_g[2,4]$ but $v_\Theta$ and $v_\vartheta$ are different there.
We take $(1/2,v_\Theta)$ as starting weight and multiplier system to define
$$A(\Gamma_g[2,4]):=\sum_{r\in\gz}[\Gamma_g[2,4],r/2,v_\Theta^r].$$
In the cases $g=1,2$ one has
$$A(\Gamma_g[2,4])=\cz[\dots f_a\dots]$$
and this ring is a polynomial ring in 2 or 4 variables. Hence the Satake compactification is $P^1$ or $P^3$.
In the case $g=3$ Runge [Ru] determined the structure of this ring.
\smallskip
To describe Runge's result,
we associate to each $a\in\fz_2^3$  a variable $F=F(a)$. We number them as follows:
\medni
\halign{\hbox to 1cm{\hfil}#)\ \hfil&(#)\quad \hfil&#)\ \hfil&(#)\quad \hfil&#) \hfil&(#)\quad \hfil&#)\ \hfil&(#)\ \hfil\cr
1&0,0,0&2&1,0,0&3&0,1,0&4&1,1,0\cr
5&0,0,1&6&1,0,1&7&0,1,1&8&1,1,1\cr}
\proclaim
{Theorem (Runge)}
{There exists a homogenous polynomial
$R\in\cz[F_1\dots F_8]$ of degree 16 such that $R$ generates the kernel of the natural homomorphism
$F_a\mapsto f_a$. This homomorphism is surjective. Hence we have
$$A(\Gamma_3[2,4]):=\sum_{r\in\gz}[\Gamma_3[2,4],r/2,v_\Theta^r]=\cz[\dots f_a\dots]=\cz[\dots F(a)\dots]/(R).$$
The subring of forms of integral weight is
$$A^{(2)}(\Gamma_3[2,4]))=\cz[\dots f_af_b\dots]=\cz[\dots\vartheta[m]^2\dots].$$
All relations between the $\vartheta[m]^2$
are  contained in the ideal generated by the Riemann relations.}
RuTh%
\finishproclaim
Runge also described the algebra of modular forms of even weight for $\Gamma_3[2]$.
We use the notation (for integral $k$ only)
$$A^{(k)}(\Gamma_g[2])=\sum_{r\equiv 0\mod k}[\Gamma_g[2],r]$$
and omit $(k)$ if $k=1$.
From Runge' s  result a modified structure theorem can be derived which we are going
to discuss now.
\smallskip
A linear subspace $G\subset \fz_2^{2g}$ is called a {\it G\"opel group\/} if it is
a maximal isotropic subspace with respect to the symplectic pairing
$$\spitz{m,n}=a'\beta+b'\alpha\qquad m=\pmatrix{a\cr b},\ n=\pmatrix{\alpha\cr\beta}.$$
One can show that each G\"opel group has a unique translate $M=a+G$ that consists
of even characteristics. The standard example is the subspace
consisting of all $\bigl({a\atop b}\bigr)$, $b=0$. Each G\"opel group
contains $2^g$ elements. We define
$$\vartheta_G=\prod_{m\in M}\vartheta[m].$$
These are modular forms on $\Gamma_g[2]$ but possibly with non-trivial multipliers.
But the following holds. Assume that $k2^{g-1}$ is divisible by 4. Then
$\vartheta_G^k$ has trivial multipliers on $\Gamma_g[2]$. We denote the smallest $k$ with this property
by $k_g$. Then  we have
$$k_g=\cases{4&if $g=1$,\cr 2&if $g=2$,\cr 1& if $g>2$.\cr}$$
More precisely we have
\proclaim
{Lemma}
{We have\smallskip
\noindent {{\rm a)}} The forms  $\vartheta_G(\tau)^{k_g}$ are in  $[\Gamma_g(2),k_g2^{g-1}]$.  \smallskip
\noindent {{b)}} The forms   $\vartheta_G(\tau)^{k_g}$ are linearly independent.}
Golm%
\finishproclaim
{\it Proof.}    The  first   statement is a consequence of the results
of  [Ig2] and it is reduced to the  fact that the  matrix  of
characteristics $M$ satisfy the congruences
$$ k_gMM'\equiv  0\,  {\rm  mod}\, 2, \qquad  {\rm diag}( k_gMM')\equiv  0\,  {\rm  mod}\, 4.$$
About the  linear independence, we  know   that the   group  $\Gamma_g$ acts
transitively   on the  even  cosets of  G\"opel groups.
Moreover, among the $\vartheta_G(\tau)^{k_g}$ only  one (namely
the standard example of all $m$ with $b=0$)
has the Fourier coefficient $a(0)\neq 0$.\qed
\medskip
In  low genera,  these  monomials are  of particular interest.
In the case $g=1$ G\"opel groups have order 2. There are three G\"opel groups which are
generated by the three elements of $\fz^2$ that are different from 0. The corresponding
$\vartheta_G$ are
$$\vartheta\Bigl[{0\atop 0}\Bigl]^4\vartheta\Bigl[{0\atop 1}\Bigl]^4,\quad
\vartheta\Bigl[{0\atop 0}\Bigl]^4\vartheta\Bigl[{1\atop 0}\Bigl]^4,
\quad \vartheta\Bigl[{0\atop 1}\Bigl]^4\vartheta\Bigl[{1\atop 0}\Bigl]^4.$$
They generate the algebra $A^{(4)}(\Gamma_1[2])$.
\smallskip
Now we consider $g=2$. There are 15 G\"opel groups. Hence we obtain 15 modular forms
$\vartheta(G)^2$ of weight $4$ one $\Gamma_2[2]$. They have trivial multipliers.
It is known that these 15 modular forms generate the ring $A^{(4)}(\Gamma_2[2])$.
The associated projective variety is the Satake compactification $\overline{\calH_2/\Gamma_2[2]}$.
It equals $\proj A^{(2)}(\Gamma_2[2])$. Due to a result
of Igusa [Ig3] this algebra is generated by $5$ forms of weight $2$ that satisfy a quartic relation.
Hence $\proj A^{(4)}(\Gamma_2[4])$ can be thought as image of the Igusa quartic in $P^4$
followed by the
Veronese embedding $P^4\to P^{14}$. Hence the variety in $ P^{14}$ can
described as the intersection of 51 quadrics,  $50$ define the Veronese embedding
and one is  induced by  Igusa's quartic.
\smallskip
Let us assume  $g=3$.
There are 135 G\"opel groups and hence 135 modular forms
$\vartheta_G$ of weight 4. They are linearly independent.
It is know through results of Runge [Ru2] that they generate the
algebra of modular forms whose weight is divisible by $4$
$$A^{(4)}(\Gamma_3[2])=\cz[\dots\vartheta_G\dots].$$
Hence $\proj(\cz[\dots\vartheta_G\dots])$ is the Satake compactification
$\overline{\calH_3/\Gamma_3[2]}$. Actally Runge proved more. The algebra
$A^{(4)}(\Gamma_3[2])$ is generated by 15 forms $\vartheta[m]^4$ and 15 of the $\vartheta_G$.
\neupara{Reciprocal maps}%
Let $\Gamma\subset\Sp(g,\gz)$ be a subgroup of finite index. Then $\calH_g/\Gamma$ carries a
structure as quasi-projective algebraic variety. Let $f_0,\dots,f_N$ be holomorphic modular
forms on $\Gamma$, different from 0,  with the same weight and the same multipliers. Then we can consider the
projective variety
$$\proj(\cz[f_0,\dots,f_N])\subset P^N.$$
We use the notation $A=\cz[f_0,\dots,f_N]$.
There is a natural rational map
$$\calH_g/\Gamma\lo\proj A,\quad Z\loma [f_0(Z),\dots,f_N(Z)].$$
It is regular outside the set of joint zeros of the forms $f_i$.
But the full (=biggest) domain of regularity can be larger of course.
Now we make the assumption that there is given a modular form $\varphi$ on $\Gamma$ (with some multiplier system)
such that
$$g_i={\varphi\over f_i}$$
is a holomorphic modular form. One can take for example the product of all $f_i$.
Then we define $B=\cz[g_1,\dots,g_N]$. This ring depends (up to canonical isomorphism)
not on the choice of $\varphi$.
We can consider
$$\calH_g/\Gamma\lo\proj B,\quad Z\loma [g_0(Z),\dots,g_M(Z)].$$
The diagram
$$\xymatrix{&\proj A\ar[dd]\\\calH_g/\Gamma\ar[ur]\ar[dr]\\& \proj B\\}$$
commutes where the vertical map is induced by the rational map
$$P^N\lo P^N,\quad [x_0,\dots,x_N]\loma [x_0^{-1},\dots,x_N^{-1}].$$
It is a birational map.
On the level of the graded algebras $A,B$ it is associated to the homomorphisms of
graded algebras
$$\eqalign{&
B\lo A,\quad g_i\loma{f_0\cdots f_N\over \varphi}g_i,\cr&
A\lo B,\quad f_i\loma {g_0\cdots g_N} f_i.\cr}$$
We want to apply this for the system of modular forms $\vartheta_G^{k_g}$. We will describe the
corresponding rings and varieties in the cases $g\le 3$.
We define
$$s_G:={\Theta_g\over \vartheta_G}\quad\hbox{where}\quad\Theta_g=\prod_m\theta[m].$$
So our pair of reciprocal rings is
$$A=\cz[\dots \vartheta_G^{k_g}\dots],\quad B=\cz[\dots s_G^{k_g}\dots].$$
\zwischen{Genus 1}%
We introduced already in the case $n=1$ for $f_i$ the three modular forms
$$f_0=\vartheta\Bigl[{0\atop 0}\Bigl]^4\vartheta\Bigl[{0\atop 1}\Bigl]^4,\quad
f_1=\vartheta\Bigl[{0\atop 0}\Bigl]^4\vartheta\Bigl[{1\atop 0}\Bigl]^4,
\quad f_2=\vartheta\Bigl[{0\atop 1}\Bigl]^4\vartheta\Bigl[{1\atop 0}\Bigl]^4.$$
They generate the algebra
$$A=A^{(4)}(\Gamma_1[2]).$$
We can take
$$\varphi=
\biggl(\vartheta\Bigl[{0\atop 0}\Bigl]\vartheta\Bigl[{0\atop 1}\Bigl]
\vartheta\Bigl[{1\atop 0}\Bigl]\biggr)^4.$$
Hence the functions $g_i$ are
$$g_0=\vartheta\Bigl[{1\atop 0}\Bigl]^4,\quad
g_1=\vartheta\Bigl[{0\atop 1}\Bigl]^4,\quad
g_2=\vartheta\Bigl[{0\atop 0}\Bigl]^4.$$
They generate the algebra
$$B=A^{(2)}(\Gamma_1[2]).$$
The relation in $A$ is
$$f_0f_1=f_1f_2+f_1f_2.$$
Hence $\proj A$ is a quadric in $P^2$ (isomorphic to $P^1$).
The relation in $B$ is $g_2=g_0+g_1$. This is a linear subspace $P^2$ (isomorphic to $P^1$).
Hence the map $\proj B\to\proj A$ is biholomorphic in this case. It is induced by the natural
embedding $A\subset B$ and it describes an isomorphism
of a quadric in $P^2$ onto a linear $P^1$.  Moreover the relations in $A$  and $B$ are reciprocally induced.
\zwischen{Genus 2}%
In the case $g=2$ we have $k_g=2$. We mentioned already
that the 15 forms $\vartheta_G^2$ generate $A^{(4)}(\Gamma_2[2])$ and that the associated
variety is a copy of the Igusa quartic. Now we study the reciprocal variety.
The $s_G$ are products of 6 theta constants. Hence $s_G^2$ are of
weight $6$. So we have 15 forms $s_G^2$ of weight 6. They are not linearly independent.
They span a space of dimension 5. There is one additional cubic relation. This defines
a copy of the Segre cubic [GS]. We obtain a birational map from the Igusa quartic and
the Segre cubic. Actually one knows that they are dual hypersurfaces. We give some details.
\proclaim
{Proposition}
{In  genus two, the   reciprocal  map induces a rational  map
$$\proj A(\Gamma_2[2]) \to P^{14}$$
which is birational unto its image.
The (closure of the) image is  defined  by
10  trinomials such as
$$s_{G_1}^2- s_{G_2}^2+s_{G_3}^2=0$$
and one cubic binomial such as
$$s_{G_1}^2 s_{G_2}^2s_{G_3}^3=s_{G_4}^2 s_{G_5}^2s_{G_6}^2.$$
The  closure of the image is isomorphic to the Segre cubic.}
G2Pr%
\finishproclaim
{\it Proof.}  This is  consequence of the results in  [GSM] and  [RSS].\qed
\smallskip
We will explain, as one can  obtain the relations  defining  the Segre cubic,
using Riemann relations.
\smallskip
We  know that, in  genus  two,  there are two types of Riemann's  quartic relations :
$$\eqalign{&
{{\rm1)}}\ \sum_{m\; \hbox{\sevenrm even}}\pm \vartheta_m^4=0,\cr
&{{\rm2)}}\ \sum_{m, n\; \hbox{\sevenrm even}}\pm \vartheta_m^2\vartheta_n^2=0.\cr}$$
The relations of the form 2)  can be obtained in this  way:  for each
one dimensional totally isotropic  subspace $N$ there are 3
even cosets $N+a=\{a, n_1+a, n_2+a, n_3+a\}$.
To each coset we associate
the  monomial
$\vartheta_{a}^2\vartheta_{a+n_1}^2$.
The  three monomials  span a two dimensional space.
Thus we  have a relation    with 3 terms. There are 15 independent such relations.
\smallskip
We shall write
$r_1+r_2=r_3$  for  such a relation.  Along $r_1r_2r_3 \neq0$ , we have
$$1/{r_2r_3}+ {1}/{r_1r_3}=  1/{r_1r_2}.$$
Multiplying by $\prod_{m,{\rm even}}\theta_m^2$, we   get a trinomial relation  as in the   proposition.
At the end we get 10 independent relations.
To  obtain the cubic relations we use   relations of the form  1). A relation looks like
$$\vartheta_{m_1}^4 \pm \vartheta_{m_2}^4 \pm\vartheta_{m_3}^4 \pm\vartheta_{m_4}^4 =0$$
with $m_1, m_2, m_3,m_4$ an azygetic  quadruplet and we apply an  argument  similar  to the previous.
We will illustrate   a specific  case.
\smallskip\noindent
We will list the  ten even characteristics.
\bigskip
\hbox{
\vbox{\halign{\qquad$#$\quad\hfill&$#$\hfill\cr
1
&(0,0,0,0)\cr
2
&(0,0,0,1)\cr
3
&(0,0,1,0)\cr
4
&(0,0,1,1)\cr
5
&(0,1,0,0)\cr}}\qquad
\vbox{\halign{$#$\quad\hfill&$#$\hfill\cr
6
&(0,1,1,0)\cr
7
&(1,0,0,0)\cr
8
&(1,0,0,1)\cr
9
&(1,1,0,0)\cr
10
&(1,1,1,1)\cr}}}
\bigni
We  have the relation
$$\vartheta_{1} ^4-\vartheta_{3}^4- \vartheta_{7}^4-\vartheta_{10}^4=0.$$
Hence
$${1}/{(\vartheta_{3} \vartheta_{7}\vartheta_{10})^4}- {1}/{(\vartheta_{ 1} \vartheta_{7}\vartheta_{10})^4}-
1/{(\vartheta_{1} \vartheta_{3}\vartheta_{10})^4}- {1}/{(\vartheta_{1} \vartheta_{3}\vartheta_{7})^4}=0.$$
Now   we multiply this relation  by
$$\prod_{i=1 }^{10} \vartheta_i^6/
((\vartheta_{ 2} \vartheta_{4}\vartheta_{5} \vartheta_{ 6})^2 \vartheta_{9}^4).$$
We observe that in  the first  denominators we  have
$$(\vartheta_{3} \vartheta_{7} \vartheta_{9}\vartheta_{10})^4(\vartheta_{ 2} \vartheta_{4}\vartheta_{5} \vartheta_{ 6})^2 =
(\vartheta_{2} \vartheta_{3}\vartheta_9\vartheta_{10})^2  (\vartheta_{4} \vartheta_{5}\vartheta_7\vartheta_{10})^2 (\vartheta_{3} \vartheta_{6}\vartheta_7\vartheta_9)^2
$$
In the  second
$$(\vartheta_{1} \vartheta_{7} \vartheta_{9}\vartheta_{10})^4(\vartheta_{ 2} \vartheta_{4}\vartheta_{5} \vartheta_{ 6})^2 =
(\vartheta_{1} \vartheta_{5}\vartheta_7\vartheta_{9})^2  (\vartheta_{1} \vartheta_{4}\vartheta_9\vartheta_{10})^2 (\vartheta_{2} \vartheta_{6}\vartheta_7\vartheta_{10})^2$$
and so on. Hence the quotients will be  of the required  form.
\smallskip
Applying  the reciprocal  map to the above relations,   one can get  the
51 quadrics  defining  $\proj A(\Gamma_2(2))^{(4)}$ in $P^{14}$.  This is  the  result of a private  communication  with  Bernd Sturmfels.\bigskip
\neupara{The G\"opel variety}%
Now we study the case $g=3$. It will be necessary to
perform some computer calculation. We used the
computer algebra system  {\saser MAGMA} [BC]. We will not
reproduce any program here, but we only give the hint
to some commands which all refer to {\saser MAGMA}.
\smallskip
Since $g=3$ we have $k_g=1$. We mentioned already that the 135
forms $\vartheta_G$ generate $A^{(4)}(\Gamma_3[2])$.
Now we want to study the reciprocal variety. For this we have to describe
the ring that is generated by the
135 forms $s_M$. They are products of 28 thetas and hence have weight 14.
We  recall that  they are  related to Coble's quartics, in fact, they
span a 15   dimensional  space that is an irreducible representation
of the group $\Gamma_3$. This  space is spanned  also  by the  15
coefficients of Coble's quartics, cf.~[GS].
\smallskip
We consider the polynomial ring in 36 variables
$$\cz[\dots T_m\dots]=\cz[T_1,\dots,T_{36}]$$
where we use the ordering of the 36 even characteristics $m$ of Sect.~1.
We also associate to each G\"opel group variables $X_G$ and $Y_G$ to build two
polynomials rings in 135 variables
$$\cz[\dots Y_G\dots],\quad \cz[\dots X_G\dots].$$
The variables $X_G$ represent the modular forms $\vartheta_G$ and the variables $Y_G$
represent the forms $s_G$. This can be expressed through the homomorphisms
$$\eqalign{&
\cz[\dots X_G\dots]\lo\cz[T_1,\dots,T_{36}],\quad X_G\loma\prod_{m\in M}T_m,\cr
&\cz[\dots Y_G\dots]\lo\cz[T_1,\dots,T_{36}],\quad Y_G\loma\prod_{m\notin M}T_m.\cr}$$
Here $M$ is the even coset of $G$.
We want to determine the inverse image of the ideal generated by the quartic Riemann relations in
$\cz[\dots Y_G\dots]$. First we describe some examples in this inverse image. We begin with linear
relations. A Riemann relation of the third kind is a relation of the form
$$r_3=r_1+r_2$$
Here the $r_i$ are products of 4 theta constants or the negative of it,
and the 12 involved characteristics are
pairwise different. Then the products $r_ir_j$ belong to the system of $X_G$ and, as a consequence,
the complementary expressions
$$Y_{ij}={\prod T_m\over r_i r_j}$$
belong to the system $Y_G$. So
$$Y_{23}=Y_{12}+Y_{13}$$
is a linear relation between the $Y_G$ which is in the inverse image of the ideal
generated by the quartic Riemann relations.
\smallskip
We construct 630 cubic relations as follows. One can find sixtuplets of G\"opel groups
$G_1,\dots,G_6$ such that
$\vartheta_{G_1}\vartheta_{G_2}\vartheta_{G_3}$ and $\vartheta_{G_4}\vartheta_{G_5}\vartheta_{G_6}$
give identical monomials in the variables $T_1,\dots,T_{36}$.
We give an example. (The digits represent theta characteristics
in the numbering defined above)
$$\eqalign{
M_1&=(3,12,13,19,24,28,30,35)\cr
M_2&=(2,8,9,12,30,32,33,36)\cr
M_3&=(5,12,16,20,21,26,30,34)\cr
\noalign{\smallskip}
M_4&=(8,12,16,19,24,26,30,33)\cr
M_5&=(3,5,9,12,30,32,34,35)\cr
M_6&=(2,12,13,20,21,28,30,36)\cr}$$
The reciprocal forms satisfy the same relation
$$s_{G_1}s_{G_2}s_{G_3}=s_{G_4}s_{G_5}s_{G_6}.$$
This gives elements of the kernel of
$\cz[\dots Y_G\dots]\lo\cz[T_1,\dots,T_{36}]$.
\smallskip
In a similar way one constructs $12\,285$ quartic relations in the kernel
of $\cz[\dots Y_G\dots]\lo\cz[T_1,\dots,T_{36}]$. There are two orbits with respect
to the full modular group.
We give two representatives.
\smallni
1) Consider the  eight G\"opel groups $G_1,\dots,G_8$ with even orbits
$$\eqalign{
M_1&=(8,12,16,19,24,26,30,33)\cr
M_2&=(6,12,13,18,24,25,31,35)\cr
M_3&=(1,11,13,19,22,27,30,35)\cr
M_4&=(1,3,6,8,25,26,27,28)\cr
\noalign{\medskip}
M_5&=(3,12,13,19,24,28,30,35)\cr
M_6&=(6,8,11,12,22,24,25,26)\cr
M_7&=(1,8,18,19,26,27,30,31)\cr
M_8&=(1,6,13,16,25,27,33,35)\cr}$$
Then we have the tautological relation
$$\vartheta_{G_1}\cdots\vartheta_{G_4}=\vartheta_{G_5}\cdots\vartheta_{G_8}.$$
The G\"opel forms $s_G$ satisfy the same relation.
\smallni
2) Consider the  eight G\"opel groups $G_1,\dots,G_8$. They agree with their even orbits, $G_i=M_i$.
$$\eqalign{
M_1&=(1,8,18,19,26,27,30,31)\cr
M_2&=(1,6,13,16,25,27,33,35)\cr
M_3&=(1,3,5,7,9,10,11,12)\cr
M_4&=1,11,15,17,24,28,32,348)\cr
\noalign{\medskip}
M_5&=(1,3,6,8,25,26,27,28)\cr
M_6&=(1,12,16,18,24,27,31,33)\cr
M_7&=(1,7,10,11,30,32,34,35)\cr
M_8&=(1,5,9,11,13,15,17,19)\cr}$$
Then we have again a tautological relation
$$\vartheta_{G_1}\cdots\vartheta_{G_4}=\vartheta_{G_5}\cdots\vartheta_{G_8}$$
and the G\"opel forms $s_G$ satisfy the same relation.
\smallskip
The space of linear relations has dimension 120.
We selected 15 independent
$Y_{G_1},\dots Y_{G_{15}}$.
The following list contains their reciprocal forms.
\medskip
\newcount\zif \zif=1
$$\halign{\qquad\the\zif#&\quad$#$\hfil\global\advance\zif by 1\cr
&T_{3}T_{12}T_{13}T_{19}T_{24}T_{28}T_{30}T_{35}\cr
&T_{2}T_{11}T_{16}T_{17}T_{23}T_{28}T_{31}T_{34}\cr
&T_{2}T_{8}T_{9}T_{12}T_{30}T_{32}T_{33}T_{36}\cr
&T_{6}T_{10}T_{13}T_{20}T_{23}T_{27}T_{32}T_{33}\cr
&T_{7}T_{10}T_{13}T_{19}T_{22}T_{27}T_{32}T_{34}\cr
&T_{5}T_{8}T_{19}T_{20}T_{21}T_{24}T_{33}T_{34}\cr
&T_{4}T_{6}T_{9}T_{12}T_{29}T_{31}T_{34}T_{35}\cr
&T_{2}T_{7}T_{17}T_{20}T_{25}T_{28}T_{30}T_{31}\cr
&T_{3}T_{4}T_{15}T_{16}T_{21}T_{22}T_{31}T_{32}\cr
&T_{3}T_{6}T_{18}T_{19}T_{25}T_{28}T_{30}T_{31}\cr
&T_{4}T_{12}T_{16}T_{17}T_{21}T_{27}T_{31}T_{34}\cr
&T_{7}T_{8}T_{15}T_{16}T_{23}T_{24}T_{29}T_{30}\cr
&T_{7}T_{9}T_{14}T_{20}T_{22}T_{28}T_{31}T_{33}\cr
&T_{4}T_{7}T_{14}T_{15}T_{26}T_{28}T_{33}T_{35}\cr
&T_{1}T_{2}T_{15}T_{16}T_{23}T_{24}T_{31}T_{32}\cr}$$
We set $Y_i=Y_{G_i}$. Then we consider the polynomial ring
$\cz[Y_1,\dots,Y_{15}]$
and the homomorphism
$$\cz[Y_1,\dots,Y_{15}]\lo\cz[\dots Y_G\dots]\lo \cz[T_1,\dots, T_{36}],\quad Y_i\loma T_{G_i}.$$
One can check that the cubic relations, considered in $\cz[Y_1,\dots,Y_{15}]$,
define a 35-dimensional space of relations. We denote the ideal generated by them in
$\cz[Y_1\cdots Y_{15}]$ by $\gota$.
The ideal generated by the cubic and quartic relations is denoted by  $\gotb$.
One can check that the ideal $\gotb$ can be generated
by the 35 cubic and by 35 quartic relations, cf [RS]. We mention that the quartic relations  can be
determined also in the following way. They have the property that
their product with $Y_1\cdots Y_{15}$ is contained in $\gota$.
\smallskip
By means of the  command {\tt quotient(a,b)}
we determined the ideal\smallskip
$$\gotx=\gota:\gotb=\{P\in\gotb;\quad P\gotb\subset\gota\}.$$
We computed that the zero locus of $\gotx$ in $P^{14}$ is zero dimensional. Hence it consists
of finitely many points in $P^{14}$. In the complement of these points the zero loci of
$\gota$ and $\gotb$ agree.
\proclaim
{Proposition}
{The variety of $\gota$ is the union of the variety of $\gotb$ and $120$ isolated points,
In particular, the variety of $\gota$ is not irreducible.}
NotIr%
\finishproclaim
The proof is given by means of a computer. \qed
\smallskip
One of the 120 points in the coordinates
$Y_1,\dots,Y_{15}$ is
$$[-1,0,1,0,1,0,1,-1,1,0,0,0,0,0,0].$$
Fortunately the variety of $\gotb$ is irreducible. Even more, $\gotb$ is a prime
ideal. The proof in [RS] seems to rely on the false assumption
that the variety of $\gota$ is irreducible. It seems us to be worthwhile to present
a correct proof. First we formulate this result.
\proclaim
{Theorem}
{The ideal  $\gotb\subset\cz[Y_1,\dots ,Y_{15}]$ is a prime ideal. It is the inverse image
of the ideal generated by the quartic Riemann relations with respect to the homomorphism
$\cz[Y_1,\dots ,Y_{15}]\to\cz[T_1,\dots,T_{36}]$. The associated
projective variety is the six-dimensional G\"opel variety.}
Goep%
\finishproclaim
{\it Proof.}
Similarly to [RS] we start with a Noether normalization which can be computed over
$\qz$ by means of the command  {\tt NoetherNormalization}.
The ideal $\gotb$ is defined over $\qz$. The ring $\qz[Y_1,\dots,Y_{15}]/\gotb_\qz$
is integral over the polynomial ring in the following 7 variables
\medni
$
3Y_{1}+2Y_{2}+Y_{3}-Y_{6}+Y_{9},\qquad
9Y_{1}+4Y_{2}+3Y_{3}+3Y_{4}+Y_{5}-Y_{8}+Y_{10},\qquad
-9Y_{1}-2Y_{2}-3Y_{3}-3Y_{4}-2Y_{5}+Y_{7}-Y_{10}+Y_{11},\qquad
11Y_{1}+5Y_{2}+4Y_{3}+4Y_{4}+Y_{5}-Y_{6}+Y_{9}+Y_{10}-Y_{11}+Y_{12},\qquad
24Y_{1}+9Y_{2}+7Y_{3}+9Y_{4}+2Y_{5}-3Y_{6}-Y_{7}-Y_{8}+2Y_{9}+2Y_{10}-2Y_{11}+Y_{12}+Y_{13},\qquad
9Y_{1}+Y_{2}+4Y_{3}+4Y_{4}+4Y_{5}+2Y_{6}-Y_{7}+Y_{8}+Y_{10}-Y_{11}+Y_{12}+Y_{14},\qquad
-25Y_{1}-8Y_{2}-6Y_{3}-9Y_{4}-Y_{5}+4Y_{6}+Y_{7}+2Y_{8}-Y_{9}-Y_{10}+2Y_{11}-Y_{13}+Y_{14}+Y_{15}.
$
\smallni
Computing Hilbert polynomials {\tt (HilbertSeries)} shows that the ring\break
$\qz[Y_1,\dots,Y_{15}]/\gotb_\qz$ is free over
this polynomial ring and hence is a Cohen-Macaulay ring (compare [RS]).
Since there can be found regular points (see below),
the ideal $\gotb$ is a radical ideal and it is equi-dimensional. It remains to show that
$\proj(\cz[Y_1,\dots,Y_{15}]/\gotb)$ has only one irreducible component. To prove this, we
intersect it with an 8-dimensional linear subspace of $P^{14}$. Actually we take the zero locus
of the last 6 linear forms in the above Noether basis. These 6 forms and $\gotb$ generate an ideal
$\gotc$. It is defined over $\qz$. One can check that it is reduced and that the dimension
of $\proj(\qz[Y_1,\dots,Y_{15}]/\gotc_\qz)$ is zero. A zero dimensional scheme is called a cluster.
The degree of this cluster is 175 (compare [RS]).
Since all irreducible components of the variety of $\gotb$ are 6-dimensional, this cluster must
hit every irreducible component. We claim now that the 175 points of this cluster are all
regular points of the variety of $\gotb$. Their explicit determination seems not to be possible,
since they have coordinates in $\bar\qz$. The situation gets better if we perform the calculation
not in $\bar\qz$ but in characteristic $p$. We took $p=557$. This reduction is allowed by a
standard flatness argument, since the
Hilbert polynomials in characteristic zero and 557 agree. The points in the algebraic closure
of $\fz_p$ can be computed by means of the command {\tt RationalPoints}. They consist of 4 Galois orbits.
Actually the orbits are
21 points in $\fz_{p^{21}}$,
22 points in $\fz_{p^{22}}$,
32 points in $\fz_{p^{32}}$,
100 points in $\fz_{p^{100}}$.
\smallskip
Since we have all points explicitly, it is no problem to show that all 175 points are regular points
of the variety of $b$. Hence each of them is contained in only one of the components of the variety
of $\gotb$. In a final step we will show that all 175 points are contained in the same component.
\smallskip
We need some more information. The forms $\vartheta_G$ and $s_G$ have level two and trivial multipliers.
Hence they must be expressible by the thetas $f_a$ of second kind. This must be reflected by homomorphims
$$\cz[\dots X_G\dots]\lo\cz[\dots F_a\dots],\quad
\cz[\dots Y_G\dots]\lo\cz[\dots F_a\dots].$$
Both homomorphisms can be constructed explicitly but the expressions of $s_G$ in terms of the $f_a$ are
extremely big (polynomials of degree 28 in 8 variables) and seem to be useless for further computations.
So we concentrate on the first of these homomorphisms and explain, how the $\vartheta_G$ can be expressed as
polynomials (of degree 8) in the $f_a$. This can be done as follows.
\smallskip
Consider a quartic Riemann relation of the kind $r_3=r_1+r_2$ where the $r_i$ are products of 4 theta constants
$\vartheta[m]$ such that the 12 involved theta constants are pairwise different. Then consider
$$2r_1r_2=r_3^2-r_1^2-r_2^2,$$
It turns out that $r_1r_2$ is one of the $\vartheta_G$. So this can be expressed by squares of theta constants,
which can be expressed by the $f_a$.
\smallskip
One can check that each of the 175 points is in general position in the sense that none of the 135
linear forms vanish on it. Hence we can invert the values of these 135 linear forms and consider
the reciprocal point
$x=[x_1,\dots,x_{135}]$.
We can consider the natural map
$$\proj(\cz[\dots f_a\dots])\lo P^{134}.$$
The image of this map is irreducible, it can be identified with
$$\proj(\cz[\dots \vartheta_G\dots])=X(\Gamma_3[2]).$$
Hence it is sufficient to prove that $x$ has an inverse image in
$$\proj(\cz[\dots f_a\dots])=X(\Gamma_3[2,4]).$$
Since the covering degree is 64, there will be 64 inverse images.
So all what one has to prove is the following.
Take the vanishing ideal of the point $x$ and take its image under the homomorphism
$\cz[\dots X_G\dots]\lo\cz[\dots F_a\dots]$. Consider the ideal that is generated by this
image and the Runge polynomial. This defines a cluster. A computer computations shows
that this cluster is not empty (and that its degree is 64).
This finishes the proof of Theorem \Goep.\qed
\vskip1cm\noindent
{\paragratit References}%
\bigskip
\item{[BC]} Bosma,\ W., Cannon,\ J., Playoust,\ C.:
{\it The Magma algebra system. I. The user language,}
J.\ Symbolic\ Comput. {\bf 24}, 235–-265 (1997)
\medskip
\item{[Co]} Coble,A.: {\it Algebraic Geometry and Theta Functions,}
American Math. Society  (1929)
\medskip
\item{[DO]} Dolgachev,\ I., Ortland,: Point Sets in Projective Spaces and Theta Functions, Ast´erisque
165 (1988).
\medskip
\item{[GS]} Grushevsky,\ S., Salvati Manni,\ R.: {\it On Coble quartic,}
Am.\ Math.\ Jornal {\bf 137}, 765--790 (2015)\hfill\break
arXiv:1212.1895
\medskip
\item{[Ig1]} Igusa,\ J.I.:
{\it On the graded ring of theta-constants,}
Am.\ J.\ of Math. {\bf 86}, 219--264 (1964)
\medskip
\item{[Ig2]} Igusa,\ J.I.:
{\it On the graded ring of theta-constant II,}
Am.\ J.\ of Math. {\bf 88}, 221--236 (1964)
\medskip
\item{[Ig3]} Igusa,\ J.I.:
{\it On Siegel Modular Forms of Genus Two (II),}
Am.\ J.\ of Math. {\bf 86}, No.\ 2, 392--412 (1964)
\medskip
\item{[Ko ]} Kondo,\ S.:
{\it  Moduli of plane quartics, G\"opel invariants and Borcherds
products,}
Int. Math. Res. Not. IMRN, {\bf12}, 2825--2860 (2011).
\item{[RS]} Ren,\ Q., Sam,\ S.V., Schrader,\ G., Sturmfels,\ B.:
{\it The universal Kummer threefold,}
Journal of
Experimental Mathematics {\bf 22}, 327--362 (2013)
\medskip
\item{[RSS]} Ren,\ Q., Sam,\ G., Sturmfels,\ B.:
{\it Tropicalization of classical moduli spaces,}
Journal of
Math. Comput. Sci. {\bf 8},119–-145, (2014)
\medskip
\item{[Ru1]} Runge. B.: {\it On Siegel modular forms, Part I,}
J.\ reine angew.\ Math. {\bf 436}, 57--84 (1993)
\medskip
\item{[Ru2]} Runge. B.: {\it On Siegel modular forms, Part II,}
Nagoya Math.\ J. {\bf 138}, 179--197 (1995)
\bye